\documentclass[12pt,
amsmath] 
{revtex4}
\usepackage{amssymb}

\renewcommand \d{\partial}

\newcommand \bQ{{\bf Q}}
\newcommand \bA{{\bf A}}
\newcommand \bX{{\bf X}}

\newcommand \bR{{\bf R}}
\newcommand \bD{{\bf D}}
\newcommand \bZ{{\bf Z}}

\newcommand \bxi{\boldsymbol{\xi}}
\newcommand \bzeta{\boldsymbol{\zeta}}
\newcommand \btheta{\boldsymbol{\theta}}
\newcommand \bfeta{\boldsymbol{\eta}}
\newcommand \bfmu{\boldsymbol{\mu}}

\newcommand \balpha{\boldsymbol{\alpha}}
\newcommand \tr{{\, \rm tr \,}}

\begin{document}

\title{
Long-term  properties of finite-correlation time isotropic stochastic systems  }
\author{ A.S. Il'yn$^{1,2}$, A.V. Kopyev$^1$, V.A. Sirota$^{1}$, and K.P.
Zybin$^{1,2}$\thanks{%
Electronic addresses:  asil72@mail.ru, kopyev@lpi.ru, sirota@lpi.ru,
zybin@lpi.ru}}
\affiliation{$^1$ P.N.Lebedev Physical Institute of RAS, 119991, Leninskij pr.53, Moscow,
Russia \\
$^2$ National Research University Higher School of Economics, 101000,
Myasnitskaya 20, Moscow, Russia}

\begin{abstract}
We consider finite-dimensional  systems of linear stochastic differential equations ${\partial
_t}{x_k}\left( t \right) = {A_{kp}}\left( t \right){x_p}\left( t \right)$, ${\bf A}(t)$ being a
stationary continuous  statistically isotropic stochastic  process with values in real $d \times d$
matrices.
  We suppose also that the laws of $\bA(t)$ satisfy the large deviation principle. For these systems, we find
    exact expressions for the Lyapunov and generalized Lyapunov exponents
    and show that they are  determined in a precise way only by the rate function of the diagonal elements of  $\bA$.

\end{abstract}

\maketitle

\section{Introduction}

Finite-dimensional first order  systems of linear differential equations with stochastic
coefficients appear in stochastic dynamics, theory of turbulence and turbulent transport, kinematic
dynamo and many other physical, chemical and biological problems (see,e.g., \cite{FGV, Sab}).
The long-term solutions of these systems 
have been investigated for the case of Gaussian and $\delta$-correlated in time stochastic noise
\cite{Kraichnan, Kazantsev}.
In this paper, we study the case of
isotropic non-Gaussian statistics and
finite correlation time of the noise.

 Consider a system of linear equations
\begin{equation}   \label{0-1}
{\partial _t}{x_k}\left( t \right) = {A_{kp}}\left( t \right){x_p}\left( t \right) , \ \ k=1\dots d
\end{equation}
where  ${\bf A}(t)$ is a stationary continuous real $d\times d$ matrix stochastic process with
known law. It is assumed to be statistically isotropic, i.e. ${\bf OA}(t){\bf O}^{T}$
   and $\bA(t)$ are identically  distributed at any $t$ and have the same time correlations  for any ${\bf O} \in O(d)$.

The formal solution of (\ref{0-1}) can be written as
$$
{x_k}\left( t \right) = {Q_{kp}}\left( t \right){x_p}\left( 0 \right),
$$
where  the evolution matrix $\bf Q$ satisfies the equation:
\begin{equation}   \label{0-2}
{\partial _t}\bQ \left( t \right) = \bA\left( t \right)\bQ \left( t \right) , \ \ \bQ \left( 0 \right) = {\bf 1}
\end{equation}
For each continuous realization of $\bA (t)$ there exists a
 well-defined solution of (\ref{0-2}) given by the Volterra product integral \cite{Gantmacher-Volterra} (in quantum-mechanical terminology  this is also called T-exponential):
\begin{equation}   \label{Volterra}
\bQ\left( t \right) = \mathop \prod \limits_{\tau  = 0}^t \left( {{\bf 1} + \bA\left( \tau  \right)d\tau } \right)
\end{equation}
So, the solutions of (\ref{0-2}) can be interpreted as continual products of random matrices.

In many physical problems (e.g., stochastic dynamics and turbulent transport)  one is basically
interested in the infinite-time asymptotic  behavior of the norms of the basis multivectors
\cite{Ottino,Pope}:
$$
{E_1} = ||\bQ{\vec e_1}|| \ , \ \ {E_2} = || \bQ{\vec e}_1 \wedge  \bQ{\vec e}_2 || , \ \ {E_d} = || \bQ{\vec e}_1 \wedge ... \wedge \bQ{\vec e}_d ||
$$
For instance, the absolute value of a passive vector advected by a random smooth flow is $E_1$; the
density of an advected passive scalar  in presence of weak molecular diffusion is $\prod _{k=1}^d
\mbox{min} \{ E_k^{-1},1 \}$ \cite{BF,PRE2017}; the energy of magnetic field generated by a
turbulent MHD flow with strong conductivity is $\mbox{min} \{ E_1^{2}/E_2^2, E_2^{2}/E_3^2 \}$
\cite{zeld, Chertkov,PhFlu21}.

From the multiplicative  Oseledets theorem \cite{Oseledets} it follows that almost surely, there
exist the limits
\begin{equation} \label{lambda-E}
\mathop {\lim }\limits_{T \to \infty }  \left( \frac 1T  \ln E_k \right)  =  \lambda_1 + \dots + \lambda _k , \ \ k=1\dots d
\end{equation}
The constants $\lambda_k$ are called Lyapunov exponents (LE).  Complete information about asymptotic
behavior of the set $\{ E_k \}$ is contained in the generalized
Lyapunov exponents (GLE) \cite{CrisantiPaladinVulp} defined by
\begin{equation} \label{GLE-E}
w\left( {{\eta_1},...,{\eta_d}} \right) = \mathop {\lim }\limits_{T \to \infty } \frac 1T \ln { \left \langle   E_{1} ^{\eta_1-\eta_2}E_{2} ^{\eta_2-\eta_3}...E_{d-1}^{\eta_{d-1}-\eta_d} E_{d}^{\eta_d} \right \rangle }
\end{equation}
where $\langle ... \rangle$ is the average over all realizations of $\bA (t)$, $\eta_k\in
\mathbb{R}$. For example, the LE can be expressed in terms of GLE:
\begin{equation} \label{lambda}
{\lambda _k} = \frac{\partial }{{\partial {\eta_k}}}w\left( 0 \right)
\end{equation}
The calculation of the values of (\ref{lambda-E}) and (\ref{GLE-E}) is the subject of this paper.

The most part of previous investigations considered the case of Gaussian  $\bA (t)$  that is
delta-correlated in time. For this case, the exact expressions for LE were found in \cite{Le Jan,
Baxendale,Newman}. The delta-correlated process typically has nowhere continuous paths, so Eq.
(\ref{0-2}) requires a stochastic convention;  the known results for that case refer to the
Stratonovich convention \cite{Gamba, Liao}.

However, in many physical applications random processes  acting as multiplicative noise are,
generally, non-Gaussian, and have finite (non-zero) correlation time. For continuous  processes the
solutions (\ref{Volterra}) are well defined, so there is no need in stochastic conventions. It
is well known \cite{Van-Kampen, JOSS1} that for such systems, the central limit theorem 'is not valid' in the sense
that all connected correlators (cumulants) of $\bA (t)$ contribute equivalently to the LE and GLE.
So, the result may differ essentially from the Gaussian  case. 

To deal with non-zero correlation time, the renovation model has been developed by different authors~\cite{renov-model1, IS-renov}: non-stationary piecewise constant process is substituted for the stationary continious process $\bA (t)$. Alternatively, we stay in the frame of stationarity and
 use the large deviation principle; in quantum mechanics and quantum field
theory the
corresponding technics is known as the low-frequency limit. 
For the processes $A_{ij}(t)$ that satisfy the large deviations
principle,  this allows to get exact expressions for LE and GLE: it appears that they are
completely determined by the rate function of the diagonal elements $A_{kk}$.



\subsection{Iwasawa decomposition}

To get explicit expressions for the LE and GLE, one performs the Iwasawa decomposition of the evolution matrix:
$$
\bQ = {\bf R  D  Z}
$$
where $\bR$ is orthogonal matrix, $\bD$ is positive diagonal, and $\bZ$ is upper-triangular unipotent matrix:
$$
{\bf R R}^T = \hat{\bf I} \ , \ \ D_{ij} = \delta_{ij} D_i , \ D_i>0 , \ \  Z_{i<j}=0 \ , Z_{jj}=1
$$
Then
$$
E_k = D_1 \dots D_k
$$
and the LE can be expressed by
\begin{equation} \label{lambda-D}
{\lambda _k} = \mathop {\lim }\limits_{T \to \infty } \frac 1T \ln {D_k} \ , \ \ \lambda _1 \ge \dots \ge \lambda _k
\end{equation}
The expression for the GLE takes the form
\begin{equation} \label{GLE-D}
w \left( {{\eta_1},...,{\eta_d}} \right) = \mathop {\lim }\limits_{T \to \infty } \frac 1T \ln { \left \langle   D_{1} ^{\eta_1}...D_{d}^{\eta_d} \right \rangle }
\end{equation}

The statistical properties of  $\bf R$, $\bf D$ and $\bf Z$ behave differently  during the evolution (\ref{0-2}); we are interested only in the statistics of $\bf D$.

In the next Section we consider the main strategy and the application of the large deviation theory on a simple one-dimensional example. Then we proceed to the multi-dimensional case, and make use of the isotropy to calculate the GLE. In conclusion, we reformulate the results in terms of the 'effective delta-process', which is a useful tool for physical applications.

\section{One-dimensional case}
\subsection{Rate function and GLE}

Consider one-dimensional differential equation with multiplicative noise:
\begin{equation} \label{1-1}
{\partial _t}x\left( t \right) = \xi \left( t \right)x\left( t \right) \ , \ \ x(0)=1
\end{equation}
where $\xi(t)$ is a continuous stationary random process with given law. We are interested in the
moments $\langle {x^\eta} \rangle$, $\eta \in \mathbb{R}$.

For any realization of $\xi(t)$,  the solution of (\ref{1-1}) is
$$
x\left( T \right) = {e^{\mathop \smallint \limits_0^T \xi \left( t \right)dt}}
$$
We assume that $\xi(t)$ satisfies the large deviation principle \cite{Varadhan}, i.e., that the
probability density of its time average
$$
\frac 1T \int \limits_0^T \xi(t) dt = \bar{\xi}
$$
satisfies at large $T\to \infty$ the relation:
$$
\rho_{\bar{\xi} } (a)  \sim  e^{-T J(a) }
$$
where the sign $\sim$ means that there exists the limit
$$
\lim \limits_{T \to \infty} \frac 1T  \ln {\rho}_{\bar{\xi} } (a)  = -  J(a)
$$
Here $J(a)$ is the rate function (Cramer function).
Then
\begin{equation} \label{x-via-w}
\langle {x^\eta(T)} \rangle = \int e^{T\eta  a} \rho_{\bar{\xi} } (a) da
\sim e^{T \, w(\eta)}
\end{equation}
where $w(\eta)$ is the Legendre transform of the rate function,
\begin{equation} \label{sup1}
w(\eta) = \sup \limits_{a} \left( \eta a-J(a) \right)
\end{equation}
This proves the existence of the limit
\begin{equation} \label{1-2}
\lim \limits_{T \to \infty}  \frac 1T \ln \langle {x^\eta} (T) \rangle  = w(\eta) \ ,
\end{equation}
$w(\eta)$ is the GLE of the process $\xi(t)$. We note that, according to (\ref{x-via-w}), the function $w(\eta)T$ at large $T$ coincides with  the cumulant generating function of the integral $\int \limits_0^T \xi(t) dt$.

Let the cumulants (or connected correlators)
\begin{equation} \label{1-2}
\left \langle \xi \left( {{t_1}} \right)...\xi \left( {{t_n}} \right) \right \rangle _c = {W^{\left( n \right)}}\left( {{t_1} - {t_2}, \ldots ,{t_1} - {t_n}} \right)
\end{equation}
be regular fast-decaying functions (i.e., let $\int W^{(n)} dt_2...dt_n$ converge). Then the cumulant-generating functional
\begin{equation} \label{Wdef}
W\left[ {\eta \left( t \right)} \right]=  \ln \left \langle {e^{\mathop \int  dt\eta  (t) \xi (t) }} \right \rangle
\end{equation}
can be   expanded into the  infinite series in $\eta(t)$:
\begin{equation} \label{1-02}
W\left[ {\eta \left( t \right)} \right] = \mathop \sum \limits_n \frac{1}{{n!}}\mathop \int  {W^{\left( n \right)}}\left( {{t_1} - {t_2}, \ldots ,{t_1} - {t_n}} \right)\eta \left( {{t_1}} \right)...\eta \left( {{t_n}} \right)d{t_1}...d{t_n}
\end{equation}

In accordance with (\ref{Wdef}),
\begin{equation} \label{const-eta}
\langle {x^\eta}\left( T \right) \rangle =  \left \langle{e^{\eta\mathop \int \limits_0^T \xi \left( t \right)dt}}\right\rangle  = e^{W\left[ {\eta{{\rm{I}}_{\left[ {0,1} \right]}}\left( t/T \right)} \right]}
\end{equation}
where ${\rm{I}}_{\left[ {0,1} \right]}$ is the indicator function of the segment $\left[ 0,1 \right]$.
Substituting (\ref{1-02}) we find
\begin{equation} \label{1-12log}
\lim \limits_{T \to \infty} \frac 1T \ln \langle {x^\eta (T)} \rangle  = \mathop \sum \limits_n   \frac{w^{(n)}}{{n!}}\eta^n
\end{equation}
where
\begin{equation} \label{1-22}
w^{(n)} =
\mathop \int  {W^{\left( n \right)}}\left( {{t_1} - {t_2}, \ldots ,{t_1} - {t_n}} \right) d{t_2}...d{t_n}
\end{equation}
Comparing (\ref{x-via-w}) with (\ref{1-12log}) we get
\begin{equation} \label{1-21}
w\left( \eta  \right) = \mathop \sum \limits_n   \frac{w^{(n)}}{{n!}}\eta^n
\end{equation}
From (\ref{1-2}), (\ref{1-21}) we see that the  long-time asymptotics of the moments
$\langle {x^\eta} \rangle$, as well as the GLE, do not depend on the details of the cumulants  and are determined only by their integrals (\ref{1-22}).

\subsection{Low-frequency limit.}
We now consider an alternative way to find the GLE based on the Lagrangian formalism. The aim is to derive the relation between the Lagrangian density and the rate function.

So, now we start from the probability density functional 
$ \wp  \left[ \xi(t) \right] $, which 
can be presented in the form
$$
 {\wp} \left[ \xi  \right]  \sim  e^{ - {\cal S} \left[ \xi  \right]} \ , \qquad
 {\cal S} \left[ {\xi \left( t \right)} \right] = \mathop \int  dt  {\cal L} (\xi, \d_t \xi, \d_t^2 \xi,\dots)
$$
Here  
 ${\cal L}$  is the Lagrangian density.
For example, the 
case
${\cal L} = \frac{1}{2}{\left( \epsilon {{\partial _\tau }\xi } \right)^2} + \frac{1}{2}{\xi ^2}$ corresponds to Ornstein-Uhlenbeck process with the only non-zero cumulant
${W^{\left( 2  \right)}}\left( {{t_1} - {t_2}} \right) = \frac{1}{2 \epsilon}{e^{ - \left| {{t_1} - {t_2}} \right|/\epsilon}}$.
 Generally,
  ${\cal L}$ is non-local, and contains the derivatives of all orders.
The generating functional can then be written in the form of the Feynman-Kac integral
\cite{Zinn-Justin}:
\begin{equation} \label{2-01}
{e^{W\left[ {\eta \left( t \right)} \right]}} \sim  \mathop \int  \left[ d\xi \right] e^{  -\int dt \left( {{\cal L}\left( t \right) - \eta \left( t \right)\xi \left( t \right)} \right)}
\end{equation}

In the previous subsection, we have seen that the  logarithm of the moments $\langle x^m \rangle$ is determined by the functions  $\eta(t)$ that are constant in the time range $\left[ 0,T\right]$ (\ref{const-eta}).  For this reason, we are interested in the values of $W\left[ {\eta \left( t \right)} \right]$  on 'slow-changing' functions with characteristic time scale $T$: $\eta_T(t) = \eta(t/T)$  
(the low-frequency limit).
In (\ref{2-01}), we rescale the time $t=\tau T$  and change accordingly the integration variable:
$
\xi \left( t \right) = {\xi _T}\left( \tau  \right)
$.
Then from (\ref{2-01}) we obtain
\begin{equation} \label{2-1}
{e^{W\left[ \eta_T \left( t\right) \right]}}\sim \mathop \int \left[ d {\xi _T} \right]
e^{-T \int d\tau {\left( {{{\cal L}_T}\left( \tau  \right) - \eta \left( \tau  \right){\xi _T}\left( \tau  \right)} \right)}   }
\end{equation}
where
$$
{{\cal L}_T}\left( \tau  \right) = {\cal L}\left(   \xi _T , \; {\frac{1}{T}{\partial _\tau }{\xi _T}} , \;
\frac{1}{T^2}{\partial _\tau^2 }{\xi _T} , \dots   \right)
$$
As $T \to \infty$, one can substitute the rate function
\begin{equation} \label{2-12}
J(\xi_T) = {{\cal L}(\xi_T,0,0,\dots)}+C
\end{equation}
for  ${{\cal L}_T}$. Here $C$ is a normalization constant. Then we estimate the integral by means
of the saddle point method:
 $$  
\mathop \int \left[ d{\xi _T} \right]
e^{-T \int d\tau {\left( J(\xi_T (\tau)) - \eta ( \tau){\xi _T} ( \tau ) \right)}   }
\sim
e^{-T \int d\tau \left( J(\xi_s(\tau))  - \eta \left( \tau  \right)  \xi _s(\tau) \right) }
$$
where $\xi_s(\tau)$ is defined by the minimum condition
\begin{equation} \label{min-cond}
\frac{{\partial J\left( {{\xi _s}} \right)}}{{\partial \xi }} = \eta \left( \tau  \right)
\end{equation}
Eventully, we get
\begin{equation} \label{2-4}
{e^{W\left[ \eta_T \left( t \right) \right]}} =
{e^{W\left[ \eta \left( t/T \right) \right]}}\sim  e^{  T \int w( \eta (\tau) ) d\tau    }
\end{equation}
where
\begin{equation} \label{2-41}
w(\eta) = \sup \limits_\xi \left( \eta \xi - J(\xi) \right)
\end{equation}
The normalization condition $W[0]=0$ results in the claim $w(0)=0$.
Now, substituting (\ref{2-4}) in (\ref{const-eta}) we find
$$
\langle {x^\eta}\left( T \right) \rangle \sim e^{T \int \limits_0^1 w(\eta) d\tau } = e^{Tw(\eta)}
$$
Comparing this with (\ref{x-via-w}), (\ref{1-2}) we see that the function $w(\eta)$ defined in this subsection coincides with the GLE found in (\ref{sup1}) and  (\ref{1-21}).  From (\ref{sup1}) and (\ref{2-41}) it also follows that the function $J(\xi)$ defined in (\ref{2-12}) coincides with the rate function.

So, both ways to determine the rate function are equivalent. In multi-dimensional case the second way appears to be more convenient.

\section{Multi-dimensional equation}

\subsection{Equation for the Iwasawa components}

Now we return to the matrix equation (\ref{0-2}) where $A(t)$ is a stationary continuous stochastic process with regular fast-decaying connected correlations:
\begin{equation} \label{3-1}
\left \langle  A_{{i_1}{j_1}}\left( {{t_1}} \right)...A_{{i_n}{j_n}} \left( {{t_n}} \right)\right \rangle_c = W_{{i_1}{j_1}...{i_n}{j_n}}^{\left( n \right)}\left( {{t_1} - {t_2}, \ldots ,{t_1} - {t_n}} \right)
 \end{equation}
and non-local Lagrangian density:
${{\cal L}_A} \left(  \bA, {{\partial _t}\bA} , \ {{\partial^2 _t}\bA} , \ \dots \right)
$.
The well-defined solution (\ref{Volterra}) exists for any continuous realization of $\bA(t)$, but
noncommutativity makes it difficult to use: there is a $T$-exponential instead of a usual exponential,
 and it seems impossible to apply the large deviations approach for $\int \bA dt$.
 However, we will see in what follows that this is possible in the case of isotropic law of $\bA$
 at least for the important 'diagonal'
 (in the sense of Iwasawa decomposition) part of the evolution matrix.

To separate the Iwasawa components, we rewrite the equation (\ref{0-2}) in the form
$$ 
 \bA = \partial_t \bQ \; \bQ^{-1}
$$
 and substitute the Iwasawa decomposition for $\bQ$.  We obtain

 \begin{equation} \label{3-4a}
\bA = \bR \, \bX \, {\bR^{ - 1}} \ , \qquad \bX =  {\bxi  + \bzeta  + \btheta }
 \end{equation}

where
 \begin{equation} \label{xi-zeta-theta-def}
\bxi  = \left( {{\partial _t}{\bf D}} \right){{\bf D}^{ - 1}} \ , \qquad  \bzeta  =
{\bf D}\left( {{\partial _t}{\bf Z}} \right){{\bf Z}^{ - 1}}{{\bf D}^{ - 1}} \ , \qquad \btheta  = {\bR^{ - 1}}\left( {{\partial _t}\bR} \right)
 \end{equation}
 The matrices $\bxi$, $\bzeta$ and $\btheta$ are diagonal, nilpotent upper triangular and antisymmetric, respectively.
The equation (\ref{xi-zeta-theta-def}) can be rewritten as
\begin{eqnarray}
{\partial _t}\bD = \bxi \bD   \label{3-5} \\
{\partial _t}\bZ = {\bD^{- 1}}\bzeta \bD \, \bZ \label{3-6} \\
{\partial _t}\bR = \bR\btheta \label{3-7}
\end{eqnarray}
Thus, treating $\bxi$, $\bzeta$ and $\btheta$ as independent variables, one could
separate the equations for $\bD$ and for $\bR$. Moreover,  the  elements $D_i$ satisfy one-dimensional equations (\ref{3-5}) same as the equation (\ref{1-1}). So,
(\ref{GLE-D})  takes the form
\begin{equation} \label{w-xi-def}
w (\eta_1,..., \eta_d) \equiv w_{\xi} (\eta_1,..., \eta_d)  = \lim \limits_{T \to \infty} \frac 1T
\ln  \left \langle e^{\int \limits_0^T (\xi_1 \eta_1 + ... + \xi_d \eta_d) dt} \right \rangle
\end{equation}
To calculate this, it would be enough to know the
 rate function of $\xi = \{ \xi_1,...,\xi_d \}$ and make use of (\ref{2-41}).
So, in the next subsection we discuss the relation between the Lagrangian densities of $\bX$ and
$\bA$.

\subsection{Change of variables}

One can consider (\ref{3-4a}) as a functional transformation from $A$ to $X$-variables,
 \begin{equation} \label{4-2}
\bA = \bR [\bX ] \bX \bR^{-1} [ \bX ] 
 \end{equation}
 where
the dependence $\bR(\bX)$ is determined by (\ref{3-7}).
To find the probability density of $\bX(t)$, one has to calculate the Jacobian:
 \begin{equation} \label{4-6}
 {\cal J} \left[ X \right] = {\rm Det }\left( {\frac{{\delta {A_{ij}}\left( t \right)}}{{\delta {X_{kp}}\left( {t'} \right)}}} \right)
 \end{equation}
It was calculated, e.g., in \cite{JOSS1}:
 \begin{equation} \label{J-exp}
{\cal J} [X] = e^{ \displaystyle\int {\tr\left( {\bfeta _0 \bX\left( t \right)} \right)dt}}
\end{equation}
where
 \begin{equation} \label{eta0}
{\left( {{\eta _0}} \right)_{kp}} = \left( {\frac{{d + 1}}{2} - k} \right){\delta _{kp}}
\end{equation}
In the Appendix A to this paper we  derive this result taking the continuous limit of a stochastic
difference equation.

With account of (\ref{4-6}) and (\ref{4-2}), from the condition ${\cal P}(\bX) \left[ d\bX\right] ={\cal P}(\bA) \left[ d \bA \right] $ we get
the  expression for the Lagrangian density of $\bX$:
$$
{\cal L}_X =  - \tr ({ \bfeta _0}\bX) + {\cal L}_A (\bA, \partial_t \bA,\partial^2_t \bA, \dots)
$$
where $\bA$ is  a function of $\bX$ in accordance with (\ref{4-2}). From (\ref{3-7}) it follows
that for any ${\bf F}(t)$,
\begin{equation} \label{deriv-of-RFR}
%
\partial_t \left( {\bR\left( t \right) {\bf F}(t) {\bR^{ - 1}}\left( t \right)} \right) = \bR\left( t \right)\left( {{D_t}{\bf F}\left( t \right)} \right){\bR^{ - 1}}\left( t \right) \ ,
 \end{equation}
where
\begin{equation} \label{long-deriv}
D_t {\bf F} = \partial_t {\bf F} + [\btheta, {\bf F}]
 \end{equation}
where $[a,b]=ab-ba$. Substituting this for the arguments of ${\cal L}_A$,  we obtain
 \begin{equation} \label{4-8}
{\cal L}_X =  - \tr ({\bf \eta _0}\bX) + {\cal L}_A
({\bf R X R}^{-1} , \bR (D_t \bX) \bR^{-1}, \bR (D^2_t \bX)
\bR^{-1}, \dots)
 \end{equation}
This expression contains $\bR$, which is the Volterra product integral of the components of $\bX$, $\bR\left[ {\bX,t} \right] = \mathop \prod \limits_{\tau  = 0}^t \left( {{\bf 1} + \btheta \left( \tau  \right)d\tau } \right)$. However,  we will see below that in the case of isotropic processes $\bR$ vanishes.

\section{Isotropic systems}

We now make use of the claim that $\bA$ is isotropic, in particular,
  $\bA(t)$ has the same probability density as ${\bf OA}(t){\bf O}^{T}$ for any  ${\bf O} \in O(d)$.
For such processes, the Lagrangian density ${\cal L}_A$  can be presented as a sum of
different combinations of traces with arguments containing products of $A$, $A^T$, and their derivatives. As one substitutes ${\bf RXR}^{-1}$ for $\bA$, the 'plates' $\bR$,$\bR^{-1}$ in the expressions vanish;
but, in accordance with (\ref{deriv-of-RFR}) and (\ref{4-8}), all time derivatives in the expression for the Lagrangian density change to 'long derivatives' (\ref{long-deriv}),
$$
\bA \mapsto \bX , \quad  \partial _t^p \bA  \mapsto  D_t^p \bX
$$
As we proceed from ${\cal L}_X$ to the rate function as in (\ref{2-12}),  we have to set all the
time derivatives of $\bX$ zero. But the commutators stay in their places, so the rate function takes
the form
\begin{equation} \label{Rate-X}
 J_X(\bX) = {\cal L}_X(\bX, 0,0,\dots) +C =
- \tr ({ \bfeta _0}\bX) + {\cal L}_A ({\bf X } , [\btheta, \bX] , [\btheta, [\btheta, \bX]], \dots)
+C
\end{equation}
One can separate the rate function of the matrix $\bA$:
$$
J_A (\bA) = {\cal L}_A (\bA, 0,\dots)
$$
and present $J_X$ in the form:
\begin{equation} \label{RateX-sum}
 J_X(\bX) = - \tr ({ \bfeta _0}\bX) + J_A (\bX) +  \delta J (\bX) +C \ ,
\end{equation}
$$
\delta J (\bX)  = {\cal L}_A ({\bf X } , [\btheta, \bX] , [\btheta, [\btheta, \bX]], \dots)-
{\cal L}_A ({\bX },0,\dots)
$$

We define  $w_X (\bfmu)$ for the matrix process $\bX$   in the same way as in
(\ref{sup1}),(\ref{2-41}),
\begin{equation} \label{Legendre-X}
w_X (\bfmu) = \sup \limits_{\bX} \left( \tr \left(\bfmu \bX\right) - J_X(\bX) \right)=\lim \limits_{T \to \infty} \frac 1T
\ln  \left \langle e^{\int \limits_0^T \tr \left(\bfmu \bX\right) dt} \right \rangle
\end{equation}
To this purpose, we have to find the extremum point $\bX_s$ analogous to $\xi_s$ in (\ref{min-cond}):
\begin{equation} \label{min-cond-matrix}
\frac {\d J_X(\bX_s)}{\d X_{qr}} = \mu_{qr}
\end{equation}
This equations system is very complicated. However,
it can be simplified significantly as we are interested only in the GLE of the diagonal elements $X_{kk}=\xi_{k}$. From 
(\ref{w-xi-def}) 
  it follows that
\begin{equation} \label{wxi-wX}
w_{\xi}(\eta_{1}, \dots, \eta_{d}) = w_X(\bfeta)  
\end{equation}
where
$\bfeta =diag(\eta_1, \eta_2,...,\eta_d)$. So, to find $w_{\xi}$ we can restrict ourself to
the diagonal matrices $\bfmu=\bfeta$ in (\ref{min-cond-matrix}), so the extremum condition
takes the form:
\begin{equation} \label{5-12}
\frac {\d J_X(\bX_s)}{\d X_{qr}} = 0 , \  q\ne r ; \qquad   \frac {\d J_X(\bX_s)}{\d X_{qq}} = \eta_{q}, \ q=1..d
\end{equation}
This  system has the diagonal solution:
\begin{equation} \label{5-13}
\bX_s = \bxi_s
\end{equation}
where ${\bxi_s} $ satisfies the relation
\begin{equation} \label{xi-s}
\frac {\d J_A (\bxi_s)}{\d X_{qq}} - (\eta_0)_{qq} = \eta_q
\end{equation}
Indeed, we recall that $\bX = \bxi + \btheta + \bzeta $, and  $J_X (\bX)$ is a combination of traces $\tr f(\bxi, \btheta, \bzeta)$ where $f$ is some product of the matrices. Since $\bxi$ is diagonal and $\btheta, \bzeta$ have zero diagonal elements, each summand in $J_X$ contains either zero or more than one of the matrices $\btheta, \bzeta$. Taking the derivative with respect to $X_{qr} , \ q\ne r$ leaves the summands with at least one multiplier  $\btheta$ or $\bzeta$,
and subsequent setting $\bX = \bxi_s$  makes them zero; so, the first equation in (\ref{5-12}) holds automatically.
On  the other hand, each summand in $\delta J (\bX)$  contains at least one $\btheta$ as a multiplier, so $\d \delta J / \d X_{qq} =0$.

We note that $\btheta_s=0$ and thus, $\delta J (\bX_s)=0$.
From (\ref{wxi-wX}), (\ref{Legendre-X}), (\ref{RateX-sum}) it then follows
\begin{equation} \label{w-intermed}
w_{\xi}(\eta_{1},\eta_{2},..., \eta_{d})
=\ tr ((\bfeta+{\bfeta_0}) {\bxi_s})   - J_A (\bxi_s) -C
\end{equation}
The statistical isotropy of $\bA$ implies $\langle \bA_{i \ne j} \rangle =0$.
Thus, $\d J_A /\d A_{qr}(\bxi_s) =0$ for $q \ne r$, and $J_A (\bxi) $ coincides  with the rate function of the diagonal elements $\alpha_{q}=A_{qq}$:
$$
 J_{\alpha} (\alpha_{1},...,\alpha_{d}) = J_A (\balpha),
$$
$$
\balpha=diag(\alpha_{1},...,\alpha_{d})
$$
 So, (\ref{w-intermed}) proves that the GLE and, as a consequence, also the Lyapunov exponents are completely determined by the rate function of the diagonal elements of $\mathbf{A}$.

By analogy with (\ref{Legendre-X}), we define the local cumulant-generating function of these
diagonal elements:
$$
w_{\alpha} (\bfeta) = \sup \limits_{\balpha} \left( \tr\left(\bfeta \balpha\right) - J_\alpha(\alpha_{1},...,\alpha_{d}) \right)=\lim \limits_{T \to \infty} \frac 1T
\ln  \left \langle e^{\int \limits_0^T \tr\left(\bfeta \balpha\right) dt} \right \rangle
$$
It is related to the cumulants of $\alpha_{i}$ by
\begin{equation} \label{cumm}
{w_\alpha }\left( {{\eta _1},...,{\eta _d}} \right) = \mathop \sum \limits_{n = 1}^\infty
 \mathop \sum \limits_{{i_1}...{i_n}} \frac{w_{{i_1}....{i_n}}^{\left( n \right)}}{{n!}}{\eta
_{{i_1}}}... {\eta _{{i_n}}}
\end{equation}
where
$$
w_{{i_1}...{i_n}}^{\left( n \right)} = \mathop \int \nolimits d{t_2}...d{t_n}
\left \langle  \alpha_{{i_1}}\left( {{t_1}} \right)...\alpha_{{i_n}} \left( {{t_n}} \right)\right \rangle_c
$$

 From (\ref{w-intermed}) we see that $w_{\xi}(\eta_1,...,\eta_d)$ can be reduced to
  $w_\alpha (\bfeta+\bfeta_0)$; with account of the normalization condition
  $w_{\xi}(0) = \ln \langle 1 \rangle = 0$, we obtain GLE and the Lyapunov exponents:
 \begin{equation} \label{5-4}
w_{\xi}(\eta_{1},\eta_{2},..., \eta_{d}) = w_{\alpha} (\bfeta+\bfeta_0) - w_{\alpha} (\bfeta_0)
\end{equation}
\begin{equation} \label{5-5}
\lambda _k = \frac{\partial }{{\partial {\eta _k}}}{w_{\alpha} }\left( \bfeta_0  \right)
\end{equation}
$$
 {\left( {{\eta _0}} \right)_{kp}} = \left( {\frac{{d + 1}}{2} - k} \right){\delta _{kp}}
$$

\subsection{Gaussian process}
Consider now the important  particular case: let $\bA (t)$ be Gaussian continuous  process with zero mean and given second-order correlator:
$$
\langle {A_{ij}}\left( t \right) \rangle = 0 \ , \qquad  \left \langle {A_{ij}}\left( {{t_1}} \right){A_{kp}}\left( {{t_2}} \right)
\right \rangle _c = K_{(ij)(kp)} \Phi (t_1-t_2)
$$
where
\begin{equation} \label{Dijkp-isotr}
K_{(ij)(kp)} =  - a{\delta _{ij}}{\delta _{kp}} + b{\delta _{ik}}{\delta _{jp}} + c{\delta _{ip}}{\delta _{jk}}
\end{equation}
where $a,b,c$ are some constants and $\Phi(t)$ is some regular even fast decaying function, $\int \Phi(t)dt = 1$.
The form (\ref{Dijkp-isotr})  of $K_{(ij)(kp)}$ is determined by the isotropy.


From (\ref{cumm}) it follows the expression for the  cumulant generation function of the diagonal
elements of $\bA$:
$$
w_{\balpha}(\bfeta)  = -\frac a2 (\tr \bfeta)^2+\frac{ b+c}2 \tr \bfeta^2
$$
From  (\ref{5-4}) we  get
\begin{equation} \label{5-14}
w_{\xi}(\eta_{1},..., \eta_{d}) =
\left( {b + c} \right)\mathop \sum \limits_{k = 1}^d \left( {\frac{{d + 1}}{2} - k} \right){\eta_k} + \frac{1}{2}\mathop \sum \limits_{k,p = 1}^d \left( {\left( {b + c} \right){\delta _{kp}} - a} \right){\eta_k}{\eta_p} \ ,
\end{equation}
\begin{equation}\label{5-15}
 \lambda_k  = \left( {b + c} \right)\left( {d + 1} - 2k \right)
\end{equation} 
In many applications one considers traceless matrices $\tr \bA = 0$. With this additional
requirement, the coefficients $a,b,c$ are associated  by the relation $b+c-ad=0$, which can be
taken into account in (\ref{5-14}).


\section{Conclusion. Effective $\delta$-process}

So, in the paper we
consider linear stochastic equations systems (\ref{0-1}) with statistically isotropic matrix random
process $\bA(t)$ that have regular fast-decaying connected correlations (\ref{3-1}). We
 find the explicit expressions (\ref{5-4}) for the generalized Lyapunov exponents in terms of
 rate functions of the diagonal elements of $\bA$.

Now we reformulate the results in the form that is useful for physical applications.

We find that the correlations of the diagonal elements of $\bA$ contribute to GLE only via their integrals:
$$
w_{{i_1}...{i_n}}^{\left( n \right)} = \mathop \smallint \nolimits d{t_2}...d{t_n}
W_{{i_1}{i_1}...{i_n}{i_n}}^{\left( n \right)}\left( {{t_1} - {t_2},...,{t_1} - {t_n}} \right)
\qquad (no \ summation)
$$
In (\ref{5-4}), the GLE are expressed in terms of the  cumulant-generating function of the
diagonal elements of the matrix~$\bA$,
$$
{w_\alpha }\left( {{\eta _1},...,{\eta _d}} \right) = \mathop \sum \limits_{n = 1}^\infty
 \mathop \sum \limits_{{i_1}...{i_n}} \frac{w_{{i_1}....{i_n}}^{\left( n \right)}}{{n!}}{\eta
_{{i_1}}}... {\eta _{{i_n}}}
$$
So, there exists the sequence of formal random processes $\bA_{\epsilon}$ with connected correlations
$$
\frac{1}{{{\epsilon ^{n - 1}}}}W_{{i_1}{j_1}...{i_n}{j_n}}^{\left( n \right)}\left( {\frac{{{t_1} -
{t_2}}}{\epsilon },...,\frac{{{t_1} - {t_n}}}{\epsilon }} \right)
$$
which produce identical GLE.  Going to the formal limit $\epsilon \to 0$ one can define the
'effective $\delta$- process' $\bA_0$ with singular correlation functions
\begin{equation} \label{last-formula}
{\rm{\Delta }}_{ij...kp}^{\left( n \right)} = w_{ij...kp}^{\left( n \right)}\delta \left( {{t_1} - {t_2}}
\right)...\delta \left( {{t_1} - {t_n}} \right)
\end{equation}
This formal process provides the same GLE and allows to split correlations and get closed equations
 for different averages (see Appendix B). Despite its formal nature, it is a handy instrument
  for calculations  for the problems that appear in theory of turbulence, turbulent transport, kinematic
   dynamo in turbulent flows etc~\cite{APS22, PoF22, PoF20}.

In the particular case of the Gaussian isotropic  processes, the result (\ref{5-15}) coincides with
the well-known expressions \cite{Liao,Gamba} obtained for the differential
equation
 $\dot{\bQ}=d{\bf W}(t)/dt \bQ$
 in the frame of Stratonovich stochastic convention, $d\bQ (t) =  d{\bf W}(t) \circ \bQ (t)$.
 Thus, for these processes the effective $\delta$-process has not only formal sense but can also be
 expressed in terms of the Wiener process' derivative. This corresponds to the Wong-Zakai theorem
 \cite{Wong-Zakai}.

We also make some notes on the relation between our approach and the renovation model~\cite{IS-renov}. In our approach, the noise is stationary for any correlation time while in the renovation model it becomes stationary only as $\tau\to 0$. Actually, the results of both approaches coincide as $\tau\to 0$. Possibly, the non-stationarity can be taken into account in our approach by means of corrections to $w_{i_1\cdots i_n}^{(n)}$. However, this is a subject for separate issue.

It is also important to note that even in the isotropic case, one can substitute the effective
$\delta$-process for the real matrix process only when calculating the long-term asymptotics of
$E_k$ and their combinations. For the quantities that depend on the non-diagonal elements
$\theta_{ij}$ and $\zeta_{ij}$ (e.g., the coordinates $x_k$), the asymptotic behavior is determined
not only by the rate function of the matrix elements $A_{ij}$ but also by the shape of their
correlation functions. This illustrates the fact that the possibility to introduce the effective
$\delta$-process is a non-trivial feature of multi-dimensional isotropic stochastic systems with
multiplicative noise.

\acknowledgments  The authors are grateful to Professor A. V. Gurevich for his permanent attention to their work.

\section*{ Appendix A: Calculation of the Jacobian (\ref{4-6})}

The Jacobian was calculated in \cite{JOSS1} by means of operations with continuous stochastic
processes and functional integrals. The notion of Jacobian is difficult to define for a continuous
process. It is natural to define it for a discrete process and then take the continuous limit. The
result must not depend on the way of discretization. Here we calculate the Jacobian for  one
particular  discretization of the random process $\bQ$ and the corresponding stochastic difference
equation. Then we show that for another choice of the discretization  with the same continuous
limit, the result is the same. In analogous way, one can check that any difference equation with
the same continuous limit would lead to the same result.

\subsection*{Discretization {\it $\grave{a}$ la} Stratonovich}

So, again we start with the equation (\ref{0-2}),
$$
{\partial _t}\bQ \left( t \right) = \bA\left( t \right)\bQ \left( t \right) , \ \ \bQ \left( 0 \right) = {\bf 1}
$$
Although $A(t)$ is a finite - correlation time process and thus the differential  equation is well defined, we still consider its discrete analog in order to substitute ordinary derivatives for variational derivatives. Thus, we
 split T into N discrete intervals, $\Delta t$ being much smaller than the correlation time of $\bA$, and consider the discrete equation:
\begin{equation} \label{Ap-1}
 \Delta \bQ_t \equiv \bQ_{t+1} - \bQ_t = \bA_t  \frac {\bQ_t +  \bQ_{t+1} }2  \Delta t
\end{equation}
We use the Stratonovich-type discretization because to calculate the Jacobian, we need the second order accuracy in $\Delta t$, and the Stratonovich choice  allows to get this accuracy for $\Delta \bQ$ without writing the second order derivative.

Multiplying (\ref{Ap-1}) by $\bQ_t^{-1}$ and by $\bQ_{t+1}^{-1}$ and taking the sum, we get:
$$
\bA_t \Delta t =  \frac 12 \left( \bQ_{t+1}  \bQ_t^{-1} -  \bQ_t \bQ_{t+1}^{-1} \right) + O(\Delta t^3)
$$
Now we make use of the Iwasawa decomosition and substitute $\bQ_t = \bR_t \bD_t \bZ_t$, $\bQ_{t+1} = \bR_{t+1} \bD_{t+1} \bZ_{t+1}$. With unified notation $\Delta {\bf F}_t = {\bf F}_{t+1}-{\bf F}_t$ and taking into account
$\Delta \left({\bf F}^{-1}\right)_t = -{\bf F}_t^{-1} \Delta  {\bf F}_t {\bf F}_t^{-1} + {\bf F}_t^{-1}\Delta  {\bf F}_t {\bf F}_t^{-1}\Delta  {\bf F}_t {\bf F}_t^{-1} + O(\Delta t^3)$
we obtain
\begin{eqnarray} \nonumber
\bA\Delta t= \frac 12 \left( \Delta \bR_t \bD_{t+1} \bZ_{t+1} +  \bR_t \Delta\bD_{t} \bZ_{t+1} + \bR_t \bD_{t} \Delta \bZ_{t} \right) \bZ_t^{-1} \bD_t^{-1} \bR_t^{-1}-
\\ \nonumber
- \frac 12  \bR_t \bD_{t} \bZ_{1} \left( \Delta (\bZ_t^{-1}) \bD_{t+1}^{-1} \bR_{t+1}^{-1} +
 \bZ_t^{-1} \Delta (\bD_{t}^{-1}) \bR_{t+1}^{-1} +  \bZ_t^{-1} \bD_{t}^{-1} \Delta (\bR_{t}^{-1})  \right)=
\\  \label{App-A-Iw}
\begin{array}{l}  
=\bR_t \left[  \frac 12 \left( \bR_t^{-1}  \bR_{t+1}  -  \bR_{t+1}^{-1} \bR_t \right) +
\frac 12 \left(  \bR_t^{-1}  \bR_{t+1} \Delta \bD_t  \bD_t^{-1} + \Delta \bD_t  \bD_t^{-1} \bR_{t+1}^{-1} \bR_t
-(\Delta \bD_t \bD_t^{-1})^2 \right) \right.+
\\    
\left. + \frac 12 \left(  \bR_t^{-1}  \bR_{t+1} \bD_{t+1} \Delta \bZ_t \bZ_t^{-1} \bD_t^{-1}   +
\bD_{t} \Delta \bZ_t \bZ_t^{-1} \bD_{t+1}^{-1} \bR_{t+1}^{-1} \bR_t  - (\bD_{t} \Delta \bZ_t \bZ_t^{-1} \bD_{t}^{-1})^2
\right)  \right] \bR_t^{-1} +O(\Delta t^3)
\end{array}
\end{eqnarray}
We denote
\begin{equation} \label{App-theta}
\btheta_t \Delta t \equiv  \frac 12 \left( \bR_t^{-1}  \bR_{t+1}  -  \bR_{t+1}^{-1} \bR_t \right)
= \frac 12 \left( \bR_t^{-1} +\bR_{t+1}^{-1} \right) \Delta \bR_t
\end{equation}
This is an antisymmetric matrix; from (\ref{App-theta}) it follows that, with accuracy $O(\Delta t^3)$,
\begin{equation} \label{App-R-theta}
\Delta \bR_t = \frac{\bR_t+\bR_{t+1}}{2} \btheta_t \Delta t
\end{equation}
which is the discrete analog to (\ref{3-7}) in accord with the Stratonovich approach. Analogously,
we claim that
$$
\Delta \bD_t =  \bxi_t  \frac{\bD_t+\bD_{t+1}}{2}\Delta t
$$
and arrive at
$$
\bxi_t \Delta t = \Delta \bD_t \bD_t^{-1} - \frac 12 \left( \Delta \bD_t \bD_t^{-1} \right)^2 + O(\Delta t^3)
$$
Demanding
$$
\Delta \bZ_t = \left( \frac{\bD_t+\bD_{t+1}}{2} \right) ^{-1} \bzeta_t \left( \frac{\bD_t+\bD_{t+1}}{2} \right)
\left( \frac{\bZ_t+\bZ_{t+1}}{2} \right) \Delta t
$$
we also get (up to $O(\Delta t^3)$ accuracy):
$$
\begin{array}{ll}
\bzeta_t \Delta t &= (\bD + \frac{\Delta \bD}2)  \Delta \bZ \left( \frac{\bZ_t+\bZ_{t+1}}{2} \right)^{-1} (\bD + \frac{\Delta \bD}2) ^{-1} \\
&=\left( {\bf 1} + \frac{\Delta \bD_t \bD_t^{-1}}{2}\right) \bD_t \Delta \bZ_t \bZ_t^{-1}  \bD_t^{-1}
\left( {\bf 1}- \frac{\bD_t \Delta \bZ_t \bZ_t^{-1}  \bD_t^{-1} }2 \right)
\left( {\bf 1} - \frac{\Delta \bD_t \bD_t^{-1}}{2}\right) +O(\Delta t^3)
\end{array}
$$
The first summand in the round brackets in (\ref{App-A-Iw}) is $\btheta_t$; now, we note that the summand in the second brackets is
$$
\bR_t^{-1} \bR_{t+1} \bxi_t \Delta t +  \bxi_t \Delta t  \bR_t^{-1} \bR_{t+1} + O(\Delta t^3) =
({\bf 1} + \btheta_t \Delta t) \bxi_t \Delta t + \bxi_t \Delta t({\bf 1} + \btheta_t \Delta t) + O(\Delta t^3)
$$
Finally, the third bracket can be written as $\bR_t^{-1} \bR_{t+1} \bzeta_t \Delta t +  \bzeta_t \Delta t  \bR_t^{-1} \bR_{t+1} + O(\Delta t^3)$.

Summarizing, we rewrite (\ref{App-A-Iw}) as
\begin{equation} \label{App-A-X}
\bA_t  = \bR_t  \left( \bX _t + \frac 12   [ \btheta _t , \bX_t ]
\Delta t
\right) \bR_t^{-1}  + O(\Delta t^2)
\end{equation}
where
$$
\bX_t = \btheta_t + \bxi_t + \bzeta_t
$$
Eq. (\ref{App-A-X}) is the discrete analog to Eq. (\ref{3-4a}).

\subsection*{Calculation of the determinant}

So, now we have to calculate the Jacobian ${\cal J}= \left| \d A_{ij\;t}/\d X_{km\;t'} \right|$.

First, we note that from (\ref{App-A-X}) and (\ref{App-R-theta}) it follows
$$
\d A_{ij\,t}/\d X_{km\,t'} = 0 \quad {\mbox for \ any} \quad   t<t'
$$
(This is the manifestation of causality)
Thus, the matrix $(\d \bA/\d \bX)_{ij,t;km,t'}$ is a block triangular matrix, and its determinant is equal to the product of the determinants of the diagonal blocks, $t=t'$:
  \begin{equation} \label{App-J-prod}
{\cal J} = \prod \limits _t \left| \frac{\d A_{ij\,t}}{\d X_{km\,t} } \right|
  \end{equation}

Second, we note that, in accordance with the causality principle, from (\ref{App-R-theta}) it follows that the value of the rotation matrix $\bR_t$ depends only  on the 'previous-time' values  of $\btheta_{t'<t}$ and does not depend on the 'simultaneous' value $\btheta_t$,
$$
\d \bR_t / \d \btheta_{t' \ge t} =0
$$
 Thus, in (\ref{App-A-X}) the derivative must be taken only over the multiplier in the square brackets.
 We now introduce the multiindices $\alpha = \{ij\}$ and the  $d^2\times d^2$ matrix
$$
 {\cal R}_{ij,km} = R_{ik} R_{jm}
 $$
 Then (\ref{App-A-X}) can be written in the form:
$$
A_{\alpha} =  {\cal R}_{\alpha \beta}  \left( X_{\beta} +   M_{\beta}  \Delta t  \right)
$$
where
$$
 {\bf M}  =  \frac 12 [\btheta, \bX] 
$$
So,
$$
\left| \frac{\d A_{\alpha}}{\d X_{\gamma} } \right| = \left| {\cal R}_{\alpha \beta } \right| \cdot
\left|  \delta_{\beta \gamma} +   \frac{ \d M_{\beta }}{\d X_{\gamma}}  \Delta t  \right|
$$
Since the matrices $\bR$ are orthogonal, $\cal R$ is also orthogonal, i.e., ${\cal R}{\cal R}^T = {\bf 1}$. Thus, $\det {\cal R} = 1$. The second determinant,  to an accuracy of  $\sim O(\Delta t)$, can be reduced to the trace of $\d {\bf M}/\d \bX$:
 \begin{equation} \label{App-det-withoutR}
\left| \frac{\d A_{\alpha}}{\d X_{\gamma} } \right| = 1 +
\left( \frac{ \d M_{\beta }}{\d X_{\gamma}}\right) \delta_{\beta \gamma} \Delta t   + O(\Delta t^2)
 \end{equation}

Finally, we make use of the fact that only the lower triangular part of $\bX$ determines the values of $\btheta$, while  the diagonal and upper triangular components are 'responsible' for $\bxi$ and $\bzeta$, correspondingly. So,
$$
\theta_{ij} = \left\{ \begin{array}{l}
X_{ij} \ \ if \ i>j \ , \\ 0 \ \ if \ i=j \ , \\ -X_{ji} \  \ if \ i<j
\end{array} \right.
$$
In particular, $\d \theta_{ij}/ \d X_{km} = 0$ if $k \le m$.
Thus, taking the derivative of $\bf M$, we obtain
$$
\tr (\d {\bf M}/ \d \bX) = \sum \limits_{i,j} \left( \frac{ \d M_{ij }}{\d X_{ij}}\right) =
\frac 12 \sum \limits_{i,j} (\theta_{ii}-\theta_{jj}) + \sum \limits_{i>j} (X_{jj} - X_{ii})
 = \sum \limits_{j=1}^d (d-2j+1)X_{jj}
$$
Combining this with (\ref{App-det-withoutR}) and (\ref{App-J-prod}), we eventually get
$$
{\cal J} = \prod \limits _t ( 1+ \sum \limits_{j=1}^d \eta_j  X_{jj} \Delta t ) + O(\Delta t^2) =
e^{\sum \limits_t \sum \limits_{j=1}^d \eta_j  X_{jj} \Delta t}
 \ , \quad
\eta_j = \frac{(d-2j+1)}2
$$
Taking the continuous limit $\Delta t \to 0$ we arrive at the integral
$$
{\cal J} = e^{\displaystyle\int \sum \limits_{j=1}^d \eta_j  X_{jj} dt }
$$
which coincides with Eq-s (\ref{J-exp}), (\ref{eta0}).

\subsection*{Discretization  {\it  $\grave{a}$ la} Ito }
The same result can be obtained in other discretization settings; However, in general case one has to keep the terms up to the second order in $\Delta t$ in the initial difference equation for $\bQ$.  Here we consider the It$\hat{\mbox{o}}$-type discretization.

We rewrite the  initial continuous differential equation (\ref{0-2}), $\d_t \bQ = \bA \bQ$, in the integral form:
$$
\bQ(t) = \bQ_0 + \int \limits_{t_0}^t \bA \bQ dt
$$
Applying this equation to the time range from $t$ to $t+\Delta t$ and solving it by means of iterations,  after two iterations we get:
$$ 
\bQ(t+\Delta t) = \bQ (t) + \bar{\bA} \bQ (t) \Delta t + \frac 12  \bar{\bA}^2 \bQ (t) \Delta t^2
$$ 
where $\bar{\bA}$ is the time average of $\bA(t)$ over the range $\Delta t$.
Basing on this equation, we write the difference equation:
\begin{equation} \label{App-Ito}
\bQ_{t+1} = \bQ_t + \left( {\bA}_t  \Delta t + \frac 12  {\bA_t}^2  \Delta t^2 \right) \bQ_t
\end{equation}
Multiplying (\ref{App-Ito}) by $\bQ_t^{-1}$ and  making use of the Iwasawa decomposition,  we present the left-hand side in the form:
  \begin{equation} \label{App-Ito-Iw}
 \bQ_{t+1} \bQ_t^{-1} = \bR_t \left( {\bf 1} + \bR_t^{-1} \Delta \bR_t  \right)
\left( {\bf 1}+ \Delta (\bD \bZ)_t (\bD \bZ)_t^{-1} \right)  \bR_t^{-1}
\end{equation}
where $\Delta F_t \equiv F_{t+1} - F_t$.

Now we formally define $\bxi_t$ by
$$
\bxi_t \Delta t = \Delta \bD_t \bD_t^{-1} - \frac 12 (\Delta \bD_t \bD_t^{-1})^2
$$
which coincides up to the second order in $\Delta t$ with
$$
\Delta \bD_t = \bD_{t+1} - \bD_t =  ( \bxi_t \Delta t + \frac 12 \bxi_t^2 \Delta t^2) \bD_t
$$
In accordance with the chosen prescription,  this difference equation corresponds to the differential equation (\ref{3-5}). Accordingly, (\ref{3-7}) corresponds to
$$
\Delta \bR_t = \bR_{t+1} - \bR_t =  \bR_t ( \btheta_t \Delta t + \frac 12 \btheta_t^2 \Delta t^2 )
$$
What about the equation (\ref{3-6}) for the upper-triangular part, with account of (\ref{3-5}) we rewrite it in the form $\d_t(\bD \bZ) = (\bxi + \bzeta) \bD \bZ$. Then, the corresponding difference equation takes the form:
$$
\Delta (\bD \bZ)_t = \bD_{t+1} \bZ_{t+1} - \bD_t \bZ_t =  \left( (\bxi_t + \bzeta_t)\Delta t   + \frac 12 (\bxi_t+\bzeta_t)^2 \Delta t^2 \right) \bD_t \bZ_t
$$
Substituting these expressions in (\ref{App-Ito-Iw}) and keeping the terms of the order of $\Delta t^2$ we obtain:
$$
\bQ_{t+1} \bQ_t^{-1} = 1+ \bR_t  \left( (\btheta_t + \bxi_t + \bzeta_t) \Delta t + (\btheta_t^2 +
(\bxi_t + \bzeta_t)^2)\Delta t^2/2 + \btheta_t (\bxi_t + \bzeta_t) \Delta t^2 \right)
 \bR_t^{-1}
$$
Combining this with (\ref{App-Ito}), we get
$$
\bA_t = \bR_t  \left( \bX_t + \frac 12 [\btheta_t , \bX_t] 
\Delta t \right)
 \bR_t^{-1} \ , \quad
\bX_t = \bxi_t + \bzeta_t + \btheta_t
$$
This result coincides with (\ref{App-A-X}) obtained in the Stratonovich convention.
The rest of the derivation is the same as in the Stratonovich case.

\section*{Appendix B. Correlation splitting for  $\delta$-process}

In this appendix, we derive the analog of the Furutsu-Novikov formula for the 
$\delta$-processes.

The Furutsu-Novikov relation for regular processes has the form \cite{Furutsu,Novikov}:
\begin{equation} \label{AppB-0}
 \langle A_{ij}(t) g[\mathbf{A}]
\rangle = \sum_{n=0}^{\infty}\frac{1}{n!}\int_{\mathbb{R}^n}\bigl\langle   A_{ij}(t) A_{i_1
j_1}(t_1)\dots A_{i_n j_n}(t_n) \bigr\rangle_c\,\biggl\langle\frac{\delta^n\,g[\mathbf{A}]}{\delta
A_{i_1j_1}(t_1)\,\dots\,\delta A_{i_nj_n}(t_n)}\biggr\rangle\mathrm{d}t_1\dots\mathrm{d}t_n
\end{equation}
for any regular functional  $g[\bA]$. For the  $\delta$-process with correlation functions
(\ref{last-formula}),
\begin{equation} \label{last-formula-again}
 {\rm{\Delta }}_{ij...kp}^{\left( n
\right)} = w_{ij...kp}^{\left( n \right)}\delta \left( {{t_1} - {t_2}} \right)...\delta \left(
{{t_1} - {t_n}} \right) \ ,
\end{equation}
 it takes the form:
\begin{equation} \label{AppB-1}
\langle A_{ij}(t) g[\mathbf{A}] \rangle = \sum_{n=0}^{\infty}\frac{1}{n!}
w^{(n+1)}_{iji_1j_1...i_nj_n}
  \biggl\langle\frac{\delta^n\,g[\mathbf{A}]}{\delta A_{ i_1j_1}(t_1)\,\dots\,\delta
A_{i_nj_n}(t_n)}\biggr\rangle
\end{equation}
Here all the variational derivatives are taken at the same moment $t$.

However, in physical applications one often has to deal with 'causal' functionals, i.e., the
functionals that depend explicitly on time and satisty the causality principle:
$$
\frac{\delta\,g[t,\mathbf{A}]}{\delta A_{ij}(t')} =0 \ \  \mbox{if}  \ \ t'>t
$$
For these functionals,
$$
\frac{\delta^n\,g[t,\mathbf{A}]}{\delta A_{i_1j_1}(t_1)\,\dots\,\delta A_{i_nj_n}(t_n)} = {\rm
I}_{[0,\infty)}(t-t_1) ...  {\rm I}_{[0,\infty)}(t-t_n) G^{(n)}_{iji_1j_1...i_n
j_n}[t,t_1,...,t_n;\bA] \ ,
$$
and
\begin{equation}\label{AppB-2}
\langle A_{ij}(t') g[t, \mathbf{A}] \rangle =  \sum_{n=0}^{\infty}\frac 1{n!}
 {w^{(n+1)}_{i\,j\,i_1j_1\dots i_nj_n}} \langle G^{(n)}_{i_1j_1\dots i_nj_n}[t,t',...,t';\bA] \rangle \ \ ,
 \ t'<t
\end{equation}
\begin{equation}\label{AppB-3}
\langle A_{ij}(t') g[t, \mathbf{A}] \rangle =  \langle A_{ ij}(t')\rangle  \langle
g[t,\bA] \rangle \ \ , \  t'>t
\end{equation}
But (\ref{AppB-1}) is inapplicable for $t=t'$ because it contains undefined values ${\rm
I}_{[0,\infty)}(0)$. So, to calculate the simultaneous correlator, we have to return to
(\ref{AppB-0}) and consider some sequence of (formal) processes with cumulants converging to
(\ref{last-formula-again}).
Thus, we choose the sequence of cumulants
$$
\frac{w^{(n)}_{i_1j_1\dots
i_nj_n}}{n}\,\sum\limits_{k=1}^n\prod\limits_{\begin{smallmatrix}{l=1}\\l\ne k
\end{smallmatrix}}^{n}\delta_\epsilon(t_k-t_l),
$$
where $\delta_{\epsilon}(t)$ are even regular functions, $\delta_{\epsilon}(t)
  \stackrel{\epsilon \to 0}{\rightarrow}  \delta(t)$, $\int \delta_{\epsilon} dt = 1$.
 With account of
$$\int^t\mathrm{d}t_1\dots\int^t\mathrm{d}t_n\prod\limits_{l=1}^{n}\delta_\epsilon(t-t_l)=\frac{1}{2^n},$$
$$\int^t\mathrm{d}t_1\dots\int^t\mathrm{d}t_n\,\delta_\epsilon(t-t_k)
\prod\limits_{\begin{smallmatrix}{l=1}\\ l\ne k
\end{smallmatrix}}^{n}
\delta_\epsilon(t_k-t_l)=\frac{1}{n}\left(1-\frac{1}{2^n}\right),$$
we arrive at
\begin{equation}\label{FN-delta}
\langle A_{ij}(t) g[t, \mathbf{A}] \rangle = \sum_{n=0}^{\infty}\frac 1{(n+1)!}
 {w^{(n+1)}_{i\,j\,i_1j_1\dots i_nj_n}} \langle G^{(n)}_{i_1j_1\dots i_nj_n}[t,...,t;\bA] \rangle
\end{equation}
This relation was presented without derivation in \cite{Klyackin}; here we
 derive it taking the formal limit.

Comparing (\ref{FN-delta}) with (\ref{AppB-2}) we see that the correlations of 'causal' functionals
are discontinuous at $t'=t$ for $\delta$-processes, and the 'naive'  convention
$I_{[0,\infty)}(0)=1/2$ is valid only for the Gaussian processes.

\end{document}